\documentclass[10pt]{amsart}
\usepackage{url,graphicx}
\usepackage[dvips]{epsfig}
\usepackage{latexsym,color}
\usepackage{amscd,amssymb,verbatim}
\newcommand{\be}{\begin{enumerate}}
\newcommand{\ee}{\end{enumerate}}

\newcommand{\rr}{{\mathbb R}}

\newcommand{\zz}{{\mathbb Z}_2}

\newcommand{\hra}{\hookrightarrow}
\newcommand{\idm}{{\textrm{id}}}
\newcommand{\ra}{\rightarrow}
\newcommand{\nin}{\noindent}
\newcommand{\pr}{\noindent{\bf Proof. }}

\newcommand{\eps}{\varepsilon}

\newcommand{\mgnv}{MG_n^\mathfrak{v}}
\newcommand{\mggnv}{MG_{g,n}^\mathfrak{v}}

\newcommand{\mgonv}{MG_{1,n}^\mathfrak{v}}
\newcommand{\tmgn}{TM_{g,n}}
\newcommand{\tmon}{TM_{1,n}}

\newcommand{\wti}{\widetilde}

\newtheorem{thm}{Theorem}[section]
\newtheorem{df}[thm]{Definition}

\newtheorem{conj}[thm]{Conjecture}

\numberwithin{equation}{section}
\numberwithin{figure}{section}

\begin{document}

\title
{Moduli spaces of metric graphs of genus~1 
with marks on vertices}

\author{Dmitry N. Kozlov}
\address{Department of Mathematics, University of Bremen, 28334 Bremen,
Federal Republic of Germany}
\email{dfk@math.uni-bremen.de}
\thanks {This research was supported by University of Bremen, 
as part of AG CALTOP}

\date\today

\begin{abstract}
In this paper we study homotopy type of certain moduli spaces of
metric graphs. More precisely, 
we show that the spaces $\mgonv$, which parametrize the isometry
classes of metric graphs of genus $1$ with $n$ marks on vertices are
homotopy equivalent to the spaces~$\tmon$, which are the moduli spaces
of  tropical curves of genus~$1$ with $n$ marked points.

Our proof proceeds by providing a sequence of explicit homotopies,
with key role played by the so-called scanning homotopy. We conjecture
that our result generalizes to the case of arbitrary genus.
\end{abstract}

\maketitle

\section{Introduction}

The moduli spaces of metric graphs with marks on vertices $\mgnv$ were
recently defined in~\cite{tms}. The main motivation for their
introduction was that they serve as universes in which the moduli
spaces of tropical  curves of genus~$g$ with $n$ marked points
$\tmgn$ can be most conveniently defined. The latter were discovered
by Mikhalkin, see~\cite{mi}, as an~important concept in the field of
tropical geometry, see~\cite{ss}. The complementary study
in~\cite{tms} was dedicated to the topological properties of the
spaces $\tmgn$, paying special attention to the case of genus~1.

Here we take one step back into the broader framework where the study
of the moduli spaces of metric graphs with marks on vertices is of
interest in its own right. We again limit ourselves to the case of
genus~1. Our main result states that the space $\mgonv$ is homotopy
equivalent to $\tmon$, which means that in this case the topological
information is already encapsulated in the tropical case. We prove
this by describing a~sequence of explicit homotopies, with the
so-called {\it scanning homotopy} playing the central role.

\section{The moduli spaces of metric graphs with marks on vertices}

For the brevity of the presentation, we shall limit ourselves to the
descriptive definitions of the concepts which we need in order to
state and to prove our results. We refer the reader to \cite{tms} for
the formally complete definitions, which at times can be somewhat
technical.

Intuitively the concept of the metric graph is rather simple: one
takes a~usual undirected graph, loops and multiple edges are
specifically allowed, and adds lengths on all the edges. While there
is a~standard way to associate a~1-dimensional CW complex to every
graph, once the edge lengths are added, there is then a~natural way to
make this CW complex into a~metric space. The isometries of the metric
graphs are the isomorphisms of the underlying graph which in addition
preserve edge lengths. When some points are marked with labels $1$
through $n$ on the metric graph, then the isometries are required to
fix the marked points as well. This allows to define the isometry
classes as equivalence classes for the corresponding equivalence
relation.  In this paper we limit ourselves to the situation where the
marks are allowed to be placed on the graph vertices only.

Let us now fix a~positive integer $n$, and let $\mgnv$ denote the set
of all isometry classes of finite metric graphs with $n$ marks on
vertices, where the vertices may have several marks. In \cite{tms} we
described a~natural way to equip this set with topology. The idea is
that given a~graph $G$, the points in a~small open neighborhood of
$[G]$ are given by the isometry classes, which have a~graph
representative obtained by a~combination of the following deformations:
\begin{itemize}
\item changing the lengths of the edges of $G$ by a~small number,
\item expanding vertices of $G$ into trees, with all the edges of the
  tree being sufficiently short,
\item if the marks on vertices are involved, distributing the marks of
  every vertex arbitrarily on the vertices of the tree which replaces
  that vertex. 
\end{itemize}
The precise definition can be found in~\cite[Section 3]{tms}.  We let
$\mggnv$ denote the subspace of $\mgnv$ consisting of the isometry
classes of connected graphs of genus~$g$.

The spaces of special interest in tropical geometry are the moduli
spaces of  tropical curves of genus~$g$ with $n$ marked points
$\tmgn$. These are the subspaces of $\mggnv$ consisting of the
isometry classes whose representatives satisfy the additional
condition that for every vertex the sum of its valency with the
cardinality of its marking list is at least~$3$.

The first tool to simplify these spaces, while preserving the homotopy
type, introduced in~\cite{tms}, was the {\it shrinking bridges} strong
deformation retraction. The way it works is quite simple: all the
bridges in the metric graph shrink to points at the speed which is
inverse proportional to the edge lengths (here a~{\it bridge} is
an~edge whose deletion increases the number of connected components,
cf.\ \cite[p.\ 11]{Die}). One can show that passing to the isometry
classes this deformation is continuous and defines a strong
deformation retraction of the the spaces $\tmgn$. The resulting spaces
consist of graphs with no bridges, which are simpler. For example,
when the genus is~$1$ we end up simply with cycles (with marked
points).

It is easy to see that the shrinking bridges strong deformation
retraction still works in the more general setting of the spaces
$\mggnv$. Hence we can choose the spaces of the bridge-free metric
graphs as the starting point of our investigation.

\section{The scanning homotopy and the homotopy type of $\mgonv$}

As mentioned above, the main focus of this paper is to understand the
homotopy type of the topological space $\mgonv$ consisting of all
isometry classes of connected metric graphs of genus~$1$ with marks
$1$ through $n$ distributed on the vertices the graph, where every
single vertex is allowed to have multiple marks.

For an arbitrary positive integer $n$, we let $\wti X_n$ denote the
space obtained from $\mgonv$ by shrinking bridges. Let $X_n$ denote
the subspace of $\wti X_n$ consisting of the isometry classes of all
cycles of total length $1$. Clearly, varying the total length of the
cycle yields a~homeomorphism
\[\wti X_n\cong X_n\times(0,\infty).\]
In particular, the space $\wti X_n$ and $X_n$ are homotopy
equivalent. For convenience we now proceed to describe the latter
space directly.

\vskip5pt \nin {\bf The points.}  The points of $X_n$ are isometry
classes of cycles of length $1$ with $n$ marked vertices. Let us
identify these cycles with a~unit circle in the plane, and let us
always assume that the vertex whose list of marks includes $1$ is
located at $(-1,0)$. The other vertices (marked or not) can be placed
on the circle arbitrarily, and to pass to the isometry classes we need
to mod out by the reflection with respect to the $x$-axis (conjugation
$\zz$-action).

\begin{figure}[hbt]
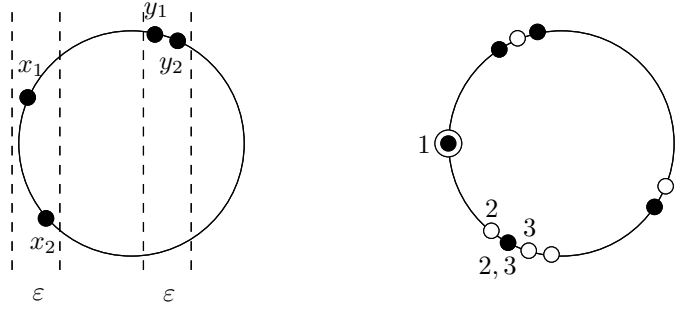

\begin{center}
  \begin{picture}(0,0)%
    \includegraphics{c1.pstex}%
  \end{picture}%
  \input{c1.pstex_t}%
  
\end{center}
\caption{On the left hand side we show pairs of $\eps$-close vertices.
On the right hand side we show two graphs $G$ and $H$ such that $H\in N_\eps(G)$.
The vertices of $G$ are filled-in and the vertices of $H$ are hollow. The vertices
which seem close on the figure are not more than $\eps$ apart.}
\label{fig:c1}
\end{figure}

\vskip5pt \nin {\bf The topology.} We say that two points
$a=(a_1,a_2)$ and $b=(b_1,b_2)$ on a~unit circle are {\it
  $\eps$-close} if the shortest path connecting $a$ and $b$ along the
circle does not leave the vertical strip
$\{(x,y)\,|\,a_1-\eps<x<a_1+\eps\}\subseteq\rr^2$. See the left hand
side of Figure~\ref{fig:c1}.

\begin{df}\label{df:top}
Let $x$ be an arbitrary point of $X_n$, and let $\eps$ be an arbitrary
positive number. We now define an~open $\eps$-neighborhood $N_\eps(x)$
of~$x$. Let $G$ be a~representative graph of~$x$. Then $y\in
N_\eps(x)$ if and only if $y$ has a~representative graph $H$ such that
\begin{enumerate}
\item[(1)] for every mark $k\in[n]$, the vertices of $G$ and of $H$
  labelled with $k$ are $\eps$-close;
\item[(2)] every vertex of $H$ is $\eps$-close to some vertex of $G$;
\item[(3)] every vertex of $G$ is $\eps$-close to some vertex of $H$. 
\end{enumerate}
\end{df}

\nin 
This description is obviously symmetric, so we have
\[x\in N_\eps(y)\Leftrightarrow y\in N_\eps(x),\]
for all $\eps>0$ and all $x,y\in X_n$. Furthermore, we set
\[N_\eps(S):=\bigcup_{x\in S}N_\eps(x),\]
for $\eps>0$ and $S\subseteq X_n$. Using these notations we have
\begin{equation}\label{eq:nee}
N_{\eps_1+\eps_2}(x)=N_{\eps_1}(N_{\eps_2}(x)),
\end{equation}
for all $\eps_1,\eps_2>0$ and $x\in X_n$.

Thinking geometrically, to obtain a~point in the $\eps$-neighborhood
of a~certain metric graph, we are allowed to shift the vertices, so
that their $x$-coordinates change by at most $\eps$, and we are
allowed to merge and to split vertices in the process. The marks
should follow with the corresponding vertices, we should merge the
mark lists when the vertices are merged, and we can split mark lists
arbitrarily when the vertices are split. See the right hand side of 
Figure~\ref{fig:c1}.

The topology on $X_n$ is now generated by the neighborhoods
$N_\eps(x)$ in the usual way: a~subspace $U$ of $X_n$ is open if and
only if for every $x\in U$ there exists $\eps>0$ so that the open
neighborhood $N_\eps(x)$ is contained in $X_n$.

Let $Y_n$ be the subspace of $X_n$ consisting of all points whose
representative graph satisfies the following conditions: 
\begin{itemize} 
\item the point with coordinates $(1,0)$ is a~vertex (which might be
  marked);
\item all other vertices of the graph are marked.
\end{itemize}
It is easy to see that the space $Y_n$ is homotopy equivalent to the
tropical moduli space $\tmon$, with the homotopy given by forgetting
the vertex at $(1,0)$, in case it is not marked.

\begin{figure}[hbt]
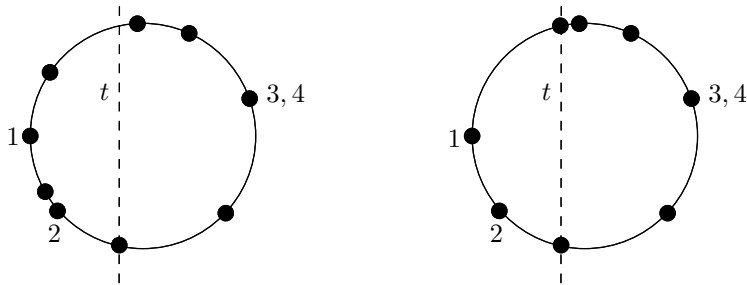

\begin{center}
  \begin{picture}(0,0)%
    \includegraphics{c2.pstex}%
  \end{picture}%
  \input{c2.pstex_t}%
  
\end{center}
\caption{On the left hand side we show a possible graph $G$; on the
  right hand side we show $\Phi(G,t)$.}
\label{fig:c2}
\end{figure}

Next, we define a~map $\Phi:X_n\times[-1,1]\ra X_n$. Let $x\in X_n$,
let $t\in[-1,1]$, and let $G$ be a~metric graph with marked vertices
representing the point~$x$. Let $H$ be the the metric graph with
marked vertices described by the following:
\begin{itemize}
\item the graph $H$ is a~cycle isometric to a~unit circle;
\item the marked vertices of $H$ are the same as those of $G$;
\item the points $(t,\sqrt{1-t^2})$ and $(t,-\sqrt{1-t^2})$ are
  vertices of $H$ (marked or not);
\item there are no unmarked vertices $(a,b)$ in $H$ satisfying $a<t$;
\item the unmarked vertices $(a,b)$ in $H$ satisfying $a>t$ are the
  same as those of the graph~$G$.
\end{itemize}
We can now set $\Phi(G,t):=H$, and accordingly $\Phi(x,t):=[H]$. It
clearly does not depend on the choice of the graph representative
of~$x$. One can visualize the homotopy defined by the map $\Phi$ as
vertical line ``scanning'' through the circle left-to-right, removing
all the unmarked vertices in the process, see Figure~\ref{fig:c2}.

\begin{thm} \label{main}
The map $\Phi$ provides a~deformation retraction from $X_n$ to $Y_n$.
In particular, the space $X_n$ is homotopy equivalent to $Y_n$, and
hence the space $\mgonv$ is homotopy equivalent to $\tmon$.
\end{thm}

\pr The main thing is to show that the set map $\Phi$ is actually
continuous.  For this purpose, choose $x\in X_n$ and $t\in[-1,1]$. Let
$\eps>0$ and choose $\alpha<\eps/2$. It would clearly to suffice to
show that for sufficiently small $\eps$ we have
\begin{equation}\label{eq:pin}
\Phi(N_\alpha(x)\times(t-\alpha,t+\alpha))\subseteq N_\eps(\Phi(x,t)).
\end{equation}

\nin
We show \eqref{eq:pin} in two steps. 

\vskip5pt

\noindent
{\bf Step 1.}  Take an arbitrary point $x\in X_n$, and let $G$ be the
metric graph representing~$x$. Furthermore, let $t_0,t_1\in[-1,1]$,
say $t_0<t_1$, such that $t_1-t_0<\eps$ for some arbitrary
$\eps>0$. What is the difference between the graphs $\Phi(G,t_0)$ and
$\Phi(G,t_1)$? The vertices, marked or not, must be the same in both
graphs, if they lie to the left of $t_0$ or the right of $t_1$, i.e.,
in the union $\{(a,b)\,|\,a<t_0\}\cup\{(a,b)\,|\,a>t_1\}$. In the
strip between $t_0$ and $t_1$ the marked vertices are the same, but
the unmarked ones may be different; see Figure~\ref{fig:c2}. However,
since the width of the strip is less than $\eps$, and since the points
of the circle on $t_0$-line are vertices of $\Phi(G,t_0)$, whereas the
points of the circle on $t_1$-line are vertices of $\Phi(G,t_1)$, we
can verify all conditions of Definition~\ref{df:top}, and conclude
that $\Phi([G],t_0)\in N_\eps(\Phi([G],t_1))$ and vice versa.

\begin{figure}[hbt]
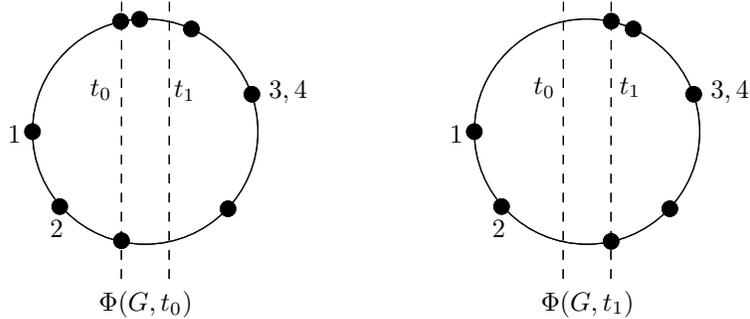

\begin{center}
  \begin{picture}(0,0)%
    \includegraphics{c3.pstex}%
  \end{picture}%
  \input{c3.pstex_t}%
  
\end{center}
\caption{We illustrate Step 1 by showing the values of $\Phi(G,-)$ for
  different times, using the graph $G$ depicted on
  Figure~\ref{fig:c2}.}
\label{fig:c3}
\end{figure}

\begin{figure}[hbt]
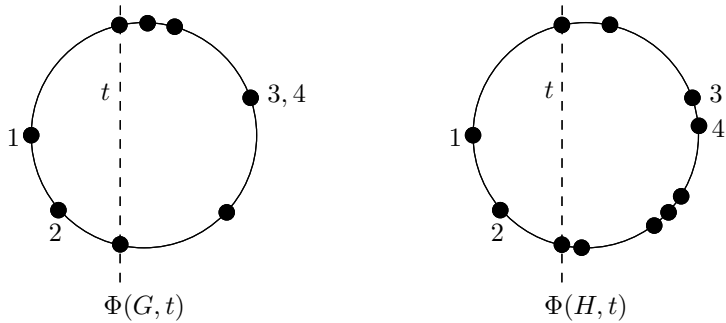

\begin{center}
  \begin{picture}(0,0)%
    \includegraphics{c4.pstex}%
  \end{picture}%
  \input{c4.pstex_t}%
  
\end{center}
\caption{We illustrate Step 2 by showing the values of $\Phi(-,t)$ for
  two close graphs.}
\label{fig:c4}
\end{figure}

\noindent
{\bf Step 2.} Fix $t\in[-1,1]$, and consider $x,y\in X_n$, and
$\eps>0$, such that $y\in N_\eps(x)$. Let $G$ be a~representative
graph for~$x$, and choose $H$ to be the representative graph for $y$
which satisfies the conditions of Definition~\ref{df:top} with respect
to the chosen graph~$G$. Again, we must ask what the difference
between the graphs $\Phi(G,t)$ and $\Phi(H,t)$ is. Specifically, we
want to prove that $\Phi([H],t)\in N_\eps(\Phi([G],t))$. To obtain the
graph $H$ from the graph $G$ we have to move vertices along the
circle, possibly merging and splitting in the process, finally getting
$\eps$-close vertices, as described in Definition~\ref{df:top}; see
Figure~\ref{fig:c4} for an~illustration.

Since the marked vertices of $\Phi(H,t)$, resp.\ $\Phi(G,t)$, are the
same as those of $H$, resp.\ $G$, the condition~(1) of
Definition~\ref{df:top} is satisfied. Furthermore, since passing from
$G$, resp.\ $H$, to $\Phi(G,t)$, resp.\ $\Phi(H,t)$, removes the
unmarked vertices to the left of the threshold line
$\{(a,b)\,|\,a=t\}$, the conditions~(2) and~(3) could theoretically be
violated, if there were points to the right of $t$ whose $\eps$-close
partner vertex has just been removed.  However, by construction the
graphs $\Phi(H,t)$ and $\Phi(G,t)$ have both points of the unit
circle, whose $x$-coordinate is $t$, as vertices. So one of these
vertices is $\eps$-close to every vertex in the strip
$\{(a,b)\,|\,t<a<t+\eps\}$. It shows, that all conditions of
Definition~\ref{df:top} are satisfied, and hence $\Phi([H],t)\in
N_\eps(\Phi([G],t))$.

We can now combine the two steps as follows. Let $x\in X_n$,
$\alpha>0$, $t\in[-1,1]$, and consider $(y,t+\delta)\in
N_\alpha(x)\times(t-\alpha,t+\alpha)$, in particular, we have
$\delta\in(-\alpha,\alpha)$. By what we proved in Step~1, we have
$\Phi(y,t+\delta)\in N_{|\delta|}(\Phi(y,t))$, while it follows from
Step~2 that $\Phi(y,t)\in N_\alpha(\Phi(x,t))$. Combined
with~\eqref{eq:nee}, these yield $\Phi(y,t+\delta)\in
N_{\alpha+|\delta|}(\Phi(x,t))$.  Since $\delta\in(-\alpha,\alpha)$
and $\alpha<\eps/2$, we conclude that $\Phi(y,t+\delta)\in
N_\eps(\Phi(x,t))$.

By construction, the map $\Phi$ provides a~homotopy between
$\Phi(-,0)=\idm_{X_n}$ and $\iota\circ\Phi(-,1):X_n\ra X_n$, where
$\iota:Y_n\hra X_n$ denotes the inclusion map. It follows that
$\Phi(-,1):X_n\ra Y_n$ is a~deformation retraction
(see~\cite[Section~6.4]{CAT}, or~\cite{Hat}), and, in particular,
$X_n$ is homotopy equivalent to~$Y_n$. As mentioned before the
theorem, the space $X_n$ is homotopy equivalent to~$\mgonv$, while the
space $Y_n$ is homotopy equivalent to $\tmon$, hence the proof is now
finished.  \qed

\vskip5pt

It is curious to note that the homotopy $\Phi$ gives a~deformation
retraction, but not a~strong deformation retraction. While being
ordinary in the classical algebraic topology, this is a~somewhat
peculiar in the context of the combinatorial algebraic topology. We
conjecture that a~stronger relation holds.

\begin{conj}\label{conj:1}
The space $Y_n$ is a~strong deformation retract of the space~$X_n$,
for all $n\geq 1$.
\end{conj}

The second conjecture is slightly more speculative, asserting that the
same holds for any genus.

\begin{conj}\label{conj:2}
The space $\tmgn$ is the strong deformation retract of $\mggnv$, for
all $g\geq 0$, and $n\geq 0$.
\end{conj}

A~natural candidate for the strong deformation retraction is provided
by the map $r:\mggnv\ra\tmgn$ which 
\begin{itemize}
\item contracts all the edges adjacent to the
  leaves\footnote{Generalizing the terminology customary for trees, we
    use the word {\it leaves} to denote any vertex of valency $1$,
    cf.\ \cite[p.\ 13]{Die}.} which are unmarked, or marked with precisely
  one label;
\item deletes all the unmarked vertices of valency~$2$.
\end{itemize}


\vskip5pt

\nin {\bf Acknowledgments.} The author is grateful Eva-Maria Feichtner
for effective discussions. He would also like to thank the anonymous
referee for the useful comments.

\end{document}